%% file: Main.tex
\documentclass[conference]{IEEEtran}
\IEEEoverridecommandlockouts
\usepackage{cite}
\usepackage{amsmath,amssymb,amsfonts,enumitem}
\usepackage{algorithmic}
\usepackage{graphicx}
\usepackage{textcomp}
\usepackage{xcolor}
\def\BibTeX{{\rm B\kern-.05em{\sc i\kern-.025em b}\kern-.08em
    T\kern-.1667em\lower.7ex\hbox{E}\kern-.125emX}}

\input{StylesMKS}

\begin{document}

\title{Assessing Power Flow Controllability via Variable Line Reactance\\
\thanks{This work was partially supported by the U.S. NSF under grant 2453151. Support from members of the Wisconsin Electric Machines and Power Electronics Consortium at the University of Wisconsin-Madison is gratefully acknowledged.}
}

\author{\IEEEauthorblockN{Eric Haag, Yuhao Chen, Giri Venkataramanan, and Manish K. Singh}
\IEEEauthorblockA{\textit{Department of Electrical \& Computer Engineering, University of Wisconsin--Madison} \\
Emails: \{erhaag, yuhao.chen,\}@wisc.edu, giri@engr.wisc.edu, manish.singh@wisc.edu}
}

\maketitle

\begin{abstract}
The rapid growth of large data center loads and inverter-based generation is increasing the stress on transmission networks, while expanding grid capacity at the required pace remains challenging. Power flow controllers (PFCs) that adjust effective line reactances to redistribute flows are often viewed as an interim solution to improve transmission network utilization. Traditional flexibility metrics and analysis approaches for PFCs focus on a limited number of operating points and contingencies. Towards gaining system-wide insights, this paper introduces a framework to quantify network flow controllability—the extent to which line flows can be reshaped through reactance adjustments. We derive analytical results demonstrating that installing PFCs on all lines enables complete controllability of feasible flow patterns. Building on these, we conduct empirical studies on the IEEE 39-bus system to examine how controllability varies with the number of PFCs and their reactance adjustment range. These analyses employ a mixed-integer linear program to optimize the siting and sizing of PFCs. Finally, we validate findings under AC power flow physics using an optimization routine that steers flows toward desired setpoints.
\end{abstract}


\section{Introduction}
The recent upsurge in large data center loads and continued increase in inverter-based resource (IBR)-connected generation have directed societal attention to the stress power grids and their operators are experiencing. While the need for increased power transmission capabilities is widely acknowledged, achieving the pace required to keep up with load growth is prohibitively challenging. Grid-enhancing technologies are gaining attention as quickly deployable interim solutions that maximize the utilization of existing transmission infrastructure~\cite{DOE2022GET,Ardakani24survey}. Broadly speaking, while system operators regulate bus power injections across a network, the line flows are determined by power flow (PF) physics. Therefore, it is frequently observed that the maximum power transfer from a source to a sink is limited by a single (or a few) lines being at their limits. Significant available power-transfer capacity on numerous power lines is rendered unusable in such cases.    

Flexible Alternating Current Transmission System (FACTS) technologies enabling control of line power flows have long been studied and developed. However, the costs and capabilities of available solutions have not been compelling enough for large-scale deployments, particularly to enable network-wide flow control. With the increase in stress on transmission infrastructure, there is a renewed interest in PF control technologies. Abstracting away implementation details, PF controllers (PFCs) on AC lines operate by varying the effective line reactance or by introducing a phase shift~\cite{Gentle2024PFC}. FACTS devices and phase-shifting transformers allow these variations, to some extent. Among these, power electronic devices that control line reactances offer the least complex approach to realize PFC, particularly in light of recent and emerging advances in power semiconductor devices with ac switching capability, converter topological features and versatility in control \cite{GaN3kV, JahnsGaN, PowerRoute, ACFacts, PowerSeControl, EvalSeries}. 

While the operation of an individual PFC is reasonably understood, system-level implications are often difficult to characterize~\cite{Tongxin16tpwrs-flexibility,Liu2018gridflexibility}. The transmission capacity enhancement benefits of PFCs are usually evaluated for specific systems, use cases, and operating scenarios~\cite{Tomsovic17tpwrs}. The currently used metrics, such as loadability margin or total transfer capability, are valuable to for providing situational insights under conceivable, actable uncertainty in potential contingencies, operating conditions, and load-change direction vectors~\cite{ChananSingh2001},~\cite{ChananSingh2002tpwrs}. However, these metrics obscure the true potential of PFCs and do not provide the systemic insights needed when uncertainty in load growth and generation installations is too high for scenario-based studies. Deviating from the conventional approaches, this work aims to directly quantify the \emph{flow controllability} offered by PFCs. We first derive analytical insights on what flow configurations can be achieved by varying line reactances within a DC PF model. After obtaining strongly positive analytical results on flow controllability enabled by installing PFCs on all lines, we proceed to empirical investigations using the IEEE 39-bus system in more practical settings. Continuing with the DC PF model, we examine the interplay among the number of PFCs installed in a network, the level of reactance control per PFC, and the achievable range of line flow alterations. These analyses were carried out using a mixed-integer linear program formulated for siting and sizing of PFCs. Finally, we validate our findings under AC PF physics using an optimal operation routine that adjusts the line reactances to steer AC power flows towards desired setpoints. Our analytical and empirical findings provide comprehensive system-level insights and set the stage for advanced analysis for optimal siting, sizing, and operation of reactance control-based PFCs.     

\emph{Notation:} Lower- (upper-) case boldface letters denote column vectors (matrices). The inequalities are understood element-wise. Vector $\be_n$ is the $n$-th canonical vector. Operators $|\cdot|$ and $\odot$ represent set cardinality and entry-wise vector products. Operator $\diag()$ yields a diagonal matrix by placing its vector argument as the main diagonal.

%%%%%%%%%%%%%%%%%%%%%%%%%%%%%%%%%%%%%%%%%%%%%%%%%%%%%%%%%
\section{Analytical Investigations}\label{sec:analytical}
In this section, we use the classical DC PF model to provide analytical insights into the controllability of power flows by varying line reactances. The formal proof of the presented results is deferred to an extended version of this manuscript due to space limitations. The single-phase equivalent of a connected power network with $N$ buses can be represented by a connected graph $\mcG=(\mcN,\mcE)$, whose nodes correspond to buses, and edges $e=(m,n)\in\mcE$ to transmission lines, with cardinality $|\mcE|=E$. We review some graph-theoretic preliminaries next that will facilitate the ensuing analysis.

\subsection{Graph-theoretic preliminaries}
 Given an undirected graph $\mcG=(\mcN,\mcE)$, one can assign arbitrary edge directions to define the $E\times N$ edge-node incidence matrix defined as
	\begin{equation}\label{eq:A}
	A_{e,k}:=
	\begin{cases}
	+1&,~k=m\\
	-1&,~k=n\\
	0&,~\text{otherwise}
	\end{cases}~\forall~e=(m,n)\in\mcE.
	\end{equation}
\begin{property}\label{prop:A1}
        For a connected network, the incidence matrix defined in~\eqref{eq:A} satisfies $\bA\bone=\bzero$ and $\rank(\bA)=N-1$.  
\end{property}

A radial topology corresponds to networks with no cycles. For such networks, $E=N-1$. Hence, the related matrix $\bA$ has full row rank per Property~\ref{prop:A1}. In a meshed graph, a sequence of adjacent edges that starts and ends at the same node is called a cycle. The number of independent cycles in a graph is given by $C=E-N+1$. For a cycle $\mcC$, one can define an indicator vector $\bn^{\mcC}\in\{0,\pm 1\}^E$ such that $n^{\mcC}_e=0$ if an edge $e$ is not included in cycle $\mcC$; $n^{\mcC}_e=1$ if edge $e$ is directed along the cycle direction, and $n^{\mcC}_e=-1$, otherwise. It holds that $\bA^\top\bn^{\mcC}=\bzero$, translating to the fact that purely cyclic network flows do not involve nodal injections or withdrawals. Let matrix $\bN\in\mathbb{R}^{E\times C}$ be constructed by stacking the indicator vectors for the $C$ independent cycles of a meshed network. It then holds that $\nullspace(\bA^\top)=\range(\bN)$.

\subsection{DC PF model-based flow controllability results}

Denote the bus active power injections by $\bp\in\mathbb{R}^N$, and line active power flows by $\bef\in\mathbb{R}^E$. 
Given the vector of line reactances, $\bx\in\mathbb{R}^E$, the DC PF equations can be stated as
\begin{subequations}\label{eq:DCPF}
    \begin{align}
        \bA^\top\bef=\bp,\label{seq:DCPFA}\\
        \diag(\bx)\bef=\bA\btheta,\label{seq:DCPFB}
    \end{align}
\end{subequations}
where $\btheta$ is the $N$-length vector of bus voltage angles. Equation~\eqref{seq:DCPFA} models the power balance at each bus, and Eq.~\eqref{seq:DCPFB} relates the voltage angle differences across a line to the power flows. Given the voltage angle at a reference bus, and power injections satisfying $\bone^\top\bp=0$, the DC power flow model is known to admit a unique solution in $(\bef, \btheta)$. However, the above is often stated in the context of networks with fixed line reactances $x_e>0$ for all edges $e$. Since we are deviating from the conventional setting of fixed positive reactances, it is worth re-examining the properties of the  DC PF model~\eqref{eq:DCPF}. The following result dictates the uniqueness of DC PF solutions when the requirement $x_e>0$ is relaxed.
\begin{proposition}\label{prop:uniqueness}
    Given network parameters $(\bA, \bx)$, balanced power injections $\bp$ such that $\bone^\top\bp=0$,  and a reference angle $\theta_r$ for some $r\in\mcN$, a unique solution $(\bef, \btheta)$ to \eqref{eq:DCPF} exists \underline{if and only if} one of the following conditions is satisfied:
    %\begin{enumerate}[label=P1.\arabic*., leftmargin=*]
    \begin{enumerate}[label=P1.\arabic*., leftmargin=3.5em]
        \item The power network is radial.
        \item For a meshed network, the matrix $\bN^\top\diag(\bx)\bN$ is invertible.
    \end{enumerate}
\end{proposition}

We will next provide results related to the following central question on power-flow controllability: 

\noindent
\textbf{Q1)} \emph{Given $(\bp,\bef)$ satisfying \eqref{seq:DCPFA}, does there always exist an $(\bx,\btheta)$ that satisfies \eqref{seq:DCPFB}?} 

Since the power injections appear only in~\eqref{seq:DCPFA}, the previous question essentially investigates the solvability of~\eqref{seq:DCPFB} written differently as
\begin{equation}\label{eq:linsys}
    \underbrace{\begin{bmatrix} \diag(\bef)& -\bA\end{bmatrix}}_{E\times (E+N)}\begin{bmatrix}  
        \bx\\ \btheta
    \end{bmatrix}=\bzero.
\end{equation}
From Property~\ref{prop:A1}, we note that $(\bx=\bzero,~\btheta=\alpha\bone)$ is a trivial solution to~\eqref{eq:linsys} for any constant scalar $\alpha$. Moreover, given that~\eqref{eq:linsys} constitutes an under-determined system of linear equations, the solution set coinciding with the nullspace of $\begin{bmatrix} \diag(\bef)& -\bA\end{bmatrix}$ has a dimension of at least $N$. Therefore, the answer to \textbf{Q1)} is trivially yes. Since transmission lines are typically inductive, a significant level of compensation would be needed to get to $x_e\leq 0$. Therefore, a pertinent follow-up question is

\noindent
\textbf{Q2)} \emph{Given line flows $\bef$, does there always exist an $\bx>\bzero$ and a vector $\btheta$ that satisfy \eqref{eq:linsys}?} 

For radial power systems, the above question readily admits a strong affirmative answer.
\begin{proposition}\label{prop:radial}
    For radial power systems, given any vector $\bef$, and any $\bx$, $\exists~\btheta$ satisfying~\eqref{eq:linsys}.
\end{proposition}
The above result follows from observing that $E=N-1$ for radial networks, resulting in linearly independent rows of $\bA$, since $\rank(\bA)=N-1$. Therefore, a solution for $\btheta$ can be obtained by solving $\bA\btheta=\diag(\bef)\bx$ for any $(\bef, \bx)$. One can therefore pick $\bx>\bzero$ to align with the pursuit of \textbf{Q2)}. However, the positive answer for radial systems has limited significance, as the notion of flow control is nuanced only for meshed networks. The next result provides a necessary and sufficient condition for answering \textbf{Q2)} with a positive. 
\begin{proposition}\label{prop:realizability}
    For meshed power systems, given any vector $\bef$, a solution $(\bx, \btheta)$ to~\eqref{eq:linsys} with $\bx>\bzero$ exists \underline{if and only if} for all $\bn\in\nullspace(\bA^\top)$ one of the following conditions hold:
    \begin{enumerate}[label=P3.\arabic*., leftmargin=3.5em]
        \item $\bef\odot\bn=\bzero$
        \item $\exists~i,j$ such that $\be_i^\top(\bef\odot\bn)>0$ and $\be_j^\top(\bef\odot\bn)<0$
    \end{enumerate}
\end{proposition}
The above conditions are not limiting, as they merely state that there cannot be cyclic flows in the network.

%%%%%%%%%%%%%%%%%%%%%%%%%%%%%%%%%%%%%%%%%%%%%%%%%%%%%%%%%%%%%%%%%%%%%%%%%
\section{Empirical Investigations}\label{sec:empirical}
The analytical results in Section~\ref{sec:analytical} establish a strong promise of redistributing power flows on transmission lines by varying the line reactances. While insightful, the findings raise further questions, such as what level of line reactance control is required to manipulate power flow patterns adequately, and how effectively power flows can be controlled when accounting for the detailed AC power flow physics. To empirically investigate these aspects, we first generated a dataset of desired line flow scenarios for the IEEE 39-bus system. Next, we formulated a DC PF-based formulation for siting and sizing of reactance control capabilities for a subset of lines such that the entire scenario set can be realized. Finally, we conducted empirical studies aiming at realizing the same set of flow scenarios via line reactance control under AC PF constraints.

%%%%%%%%%%%%%%%%%%%%%%%%%%%%%%%%%%%%%%%%
\subsection{Flow scenario set generation}
We obtained the nominal (re)active power demands ($\bp^{\mrd}_{\circ},\bq^{\mrd}_{\circ})$ for the IEEE 39-bus system~\cite{MATPOWER}. Using entry-wise scaling $\bp^\mrd_s=\bp^\mrd_{\circ}\odot\bu_s^{\mrp}$ and $\bq^\mrd_s=\bq^{\mrd}_{\circ}\odot\bu_s^{\mrq}$, $S$ demand scenarios were generated. The entries of $\bu_s^{\mrp}$ and $\bu_s^{\mrq}$ were each independently drawn from uniform random distributions on the intervals $[0.2,2]$ and $[0.9,1.1]$, respectively. To mimic a setup in which approximate grid dispatch could be carried out by relaxing the DC PF physics~\eqref{seq:DCPFB} and merely respecting resource and network-flow constraints~\eqref{seq:DCPFA}, we solve the following optimal power flow problem for all $S$ demand instances.
\begin{subequations}\label{eq:P1}
\begin{align}
\min\qquad (\bp^{\mrg}_s)^\top &\bQ\bp^\mrg_s\quad\quad\qquad\quad \quad\quad\textrm{(P1)}\notag\\
\textrm{s.to}~\bA^\top\bef_s&=\bp^\mrg_s-\bp^\mrd_s,\\
\underline{\bp}^\mrg&\leq\bp^\mrg_s\leq\bar{\bp}^\mrg,\\
-\bar{\bef}&\leq\bef_s\leq\bar{\bef},
\end{align}
\end{subequations}
where the cost coefficients $\bQ$, generation limits $(\underline{\bp}^\mrg, \bar{\bp}^\mrg)$, and line limits $\bar{\bef}$ were obtained from the MATPOWER case file~\cite{MATPOWER}. After disregarding infeasible instances, we remove scenarios that violate the conditions identified in Proposition~\ref{prop:realizability} since these instances would not be realizable with any positive value of line reactances. Thus obtained dataset $\{\bp^\mrd_s,\bq^\mrd_s,\bp^\mrg_s,\bef_s\}_{s=1}^S$ will be used in the ensuing analysis with a special focus on flows $\bef_s$. The reasoning is that the flows obtained from the above procedure are representative of desirable flow patterns on the grid. However, our analysis approach will apply to an arbitrary set of desired flow configurations.   

We next assessed how distinct the desired flow configurations $\{\bef_s\}$ are from the flows expected from DC PF. Specifically, we solved~\eqref{eq:DCPF} for all the scenarios using the nominal line reactances, denoted as $\bx^\circ$, and setting the net power injection $\bp$ as $\bp^\mrg_s-\bp^\mrd_s$ to obtain the DC PF solutions $\{\bef^{\circ}_s\}_{s=1}^S$. We then sorted the scenarios in increasing order of $\|\bef_s-\bef^\circ_s\|_2^2$, providing the intuition that, as one parses through the dataset, they encounter flow configurations that become increasingly difficult to realize by varying the line reactances.

%%%%%%%%%%%%%%%%%%%%%%%%%%%%%%%%%%%%%%%%
\subsection{DC PF based siting and sizing of reactance control}
Let the maximum and minimum possible reactance values for all lines be stacked in vectors $\bar{\bx}$ and $\underline{\bx}$, respectively. Let $\bd\in\{0,1\}^E$ capture reactance controller placement such that $d_e=1$ indicates presence of some level of reactance control while $d_e=0$ implies no reactance control, i.e., $x^{\circ}_e=\bar{x}_e=\underline{x}_e$. The optimal siting and sizing task for reactance controllers to ensure the realizability of flows $\{\bef_s\}_{s=1}^{S'}$ can be formulated as
\begin{subequations}\label{eq:P2}
\begin{align}
\min \bone^\top(\bar{\bx}-\underline{\bx})&\qquad\qquad\quad\quad\quad\quad\textrm{(P2)}\notag\\
\textrm{over}~\{\bx_s, \btheta_s&\}_{s=1}^{S'}, \bd, \underline{\bx}, \bar{\bx}\notag\\ 
\textrm{s.to}~\diag(\bef_s)\bx_s&=\bA\btheta_s,~\forall~s\label{seq:P2a}\\
\underline{\bx}&\leq\bx_s\leq\bar{\bx},~\forall~s\label{seq:P2b}\\
\mathbf{0}&\leq\bar{\bx}-\bx^\circ\leq \bx^\circ\odot\bd,\label{seq:P2c}\\
\mathbf{0}&\leq\bx^\circ-\underline{\bx}\leq \bx^\circ\odot\bd,\label{seq:P2d}\\
\bone^\top\bd&\leq K,~
\bd \in \{0,1\}^{E}\label{seq:P2f},
\end{align}
\end{subequations}
where $S'\leq S$ is the number of scenarios considered, $K$ is the budget for the number of reactance controllers. The upper bounds in~\eqref{seq:P2c}--\eqref{seq:P2d} limit the reactance control range to $\pm 100\%$, eliminating the possibility of negative line reactances. 

\begin{figure}[h]
    \centering
    \includegraphics[width=0.9\linewidth]{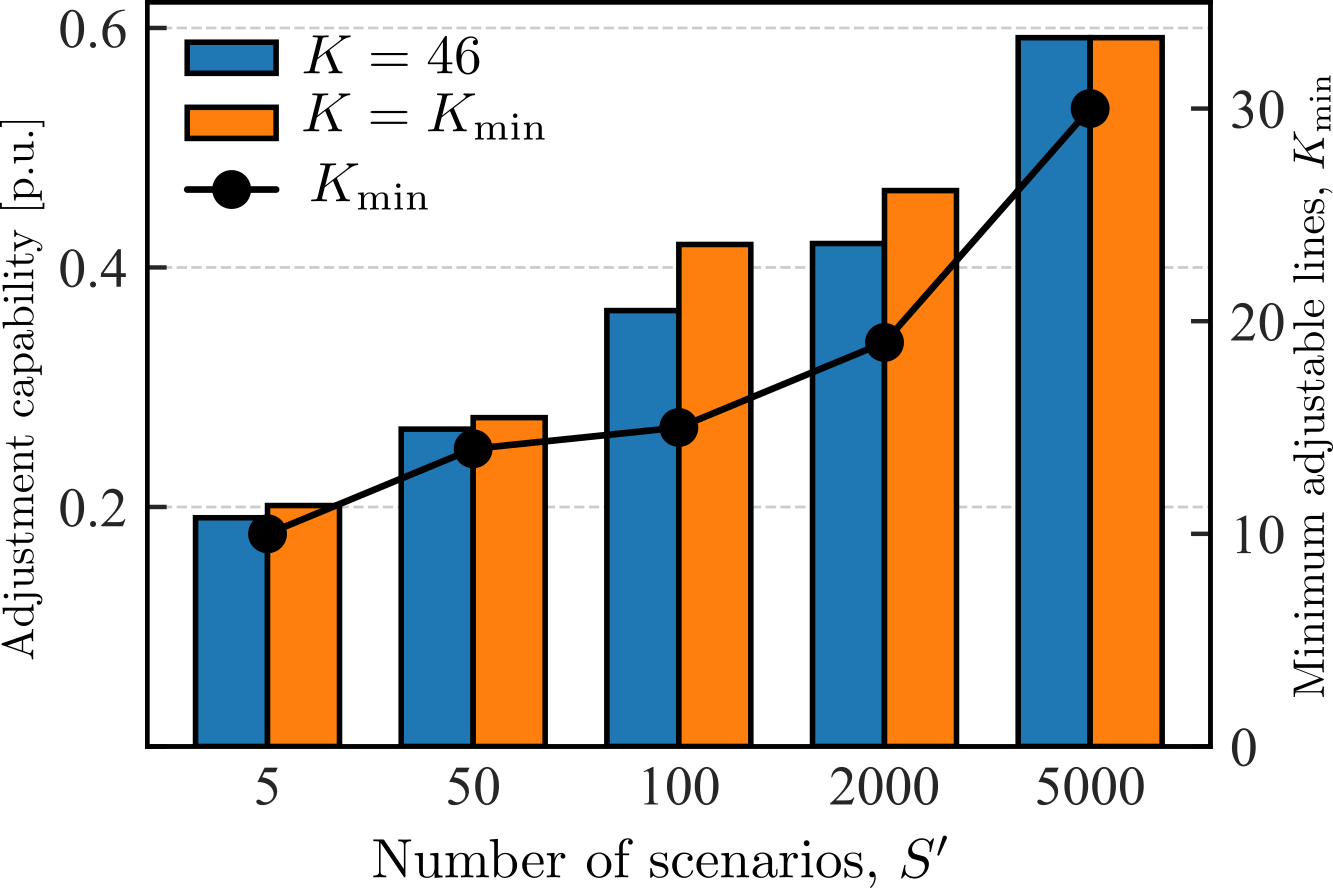}
    \caption{The trendline shows the minimum number of lines with adjustable reactance needed to realize $S'$ flow scenarios. The bar graphs provide total reactance-adjustment capacity $\bone^\top(\bar{\bx}-\underline{\bx})$ required when all 46 lines or just $K_{\min}$ lines have reactance control capability. }
    \label{fig: Kmin and adjustment needed}
\end{figure}

Using (P2), we studied the interplay between the difficulty of flow scenarios considered and the number and capacity of reactance controllers needed. Since the dataset is sorted in increasing order of difficulty surrogate $\|\bef_s-\bef^\circ_s\|_2^2$, we solved (P2) for $S'=5, 50, 100, 2000,$ and $5000$ scenarios. For each $S'$, we first solved (P2) with $K=46$, the total number of lines in the IEEE 39-bus system, and recorded the optimal total control capability requirement $\bone^\top(\bar{\bx}-\underline{\bx})$ in per unit. Then we re-solved (P2) for decreasing values of $K$ until the problem becomes infeasible. This analysis provides the minimum number of controllers needed for each $S'$, denoted as $K_{\min}(S')$, and the corresponding value of $\bone^\top(\bar{\bx}-\underline{\bx})$. Figure~\ref{fig: Kmin and adjustment needed} summarizes our findings for the above analysis.

Aligning with expectations, Fig.~\ref{fig: Kmin and adjustment needed} indicates that the minimum number of reactance controllers required increases with the number and difficulty of flow patterns to be realized. Interestingly, all flow patterns could be realized with reactance controllers placed on $30$ out of $46$ lines. Further analytical investigation revealed that reactance controllers are never needed on lines that are not in any cycles, referred to as bridges in graph terminology; a result that follows from extension of Proposition~\ref{prop:radial}. The test system had $11$ bridges. Thus, the generated flow patterns were rich enough to require controllers on $30$ out of $35$ non-bridge lines. When evaluating the total reactance control capability across all lines, Fig.~\ref{fig: Kmin and adjustment needed} shows the expected increase with the increasing number of scenarios. However, a potentially surprising observation was made when comparing the total control capability required with $K=46$ versus $K=K_{\min}$. Specifically, for any value of $S'$, the total control capability $\bone^\top(\bar{\bx}-\underline{\bx})$ does not increase significantly when one decreases the allowable number of controllers. Further analytical and empirical analysis is needed to uncover the underlying reasons for this observation.  

\begin{figure}[t]
    \centering
    \includegraphics[width=0.8\linewidth]{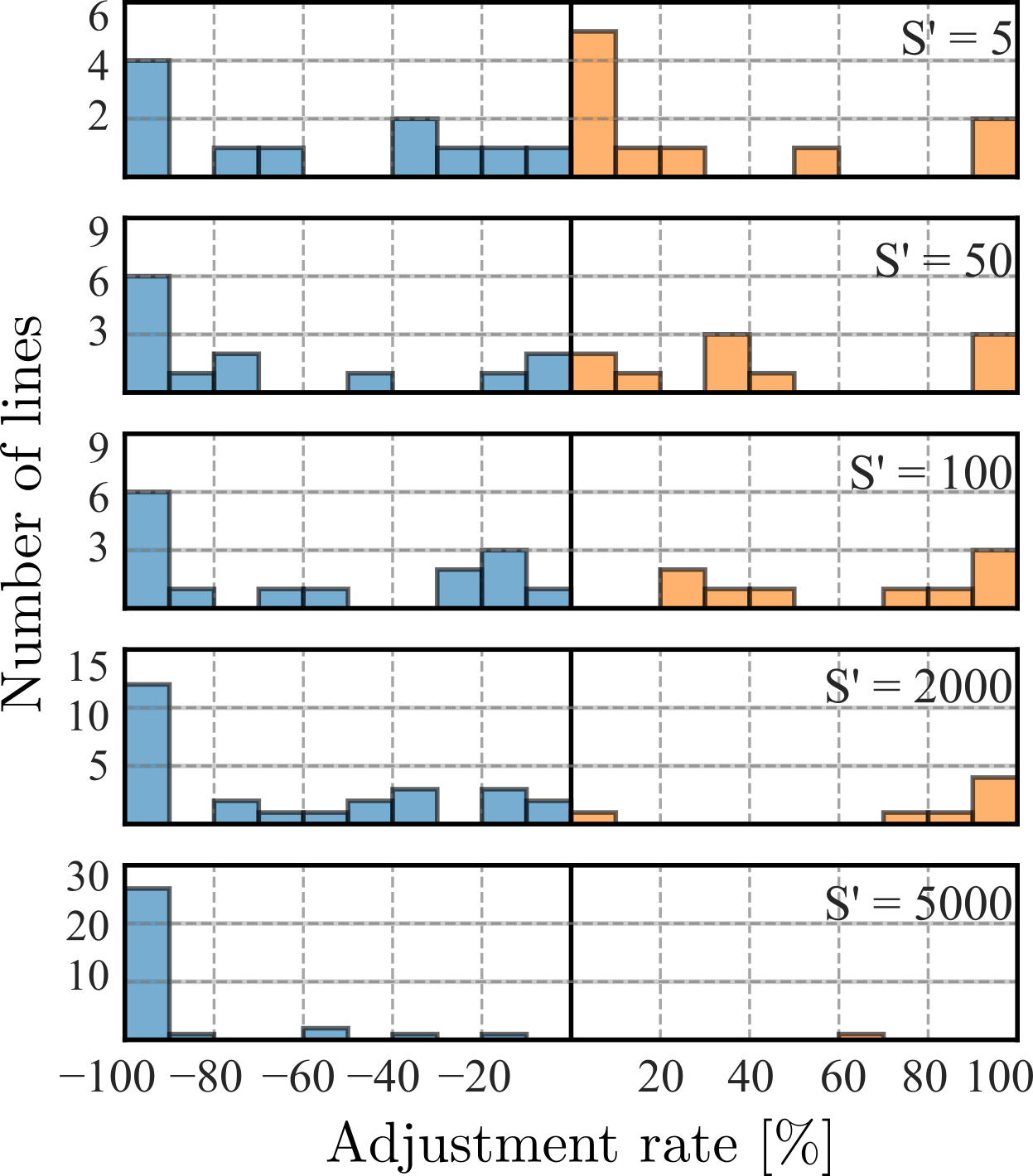}
    \vspace{-8pt}   
    \caption{Distributions of negative ($\gamma_{e,\mathrm{down}}$) and positive 
($\gamma_{e,\mathrm{up}}$) reactance-adjustment capacity needed per line for different numbers of scenarios $S'$. Despite allowing $K=46$, several lines did not need reactance adjustment. Those lines are omitted from this figure.}
    \label{fig: individual adjustment}
\end{figure}
We next studied the optimally allocated reactance control capability on individual lines for $K=46$ and while varying $S'$ using the relative reactance adjustment metrics 
\begin{equation*}
    \gamma_{e,\text{down}}=\frac{\underline{x}_e-x_{e}^\circ}{x_{e}^\circ},\qquad\gamma_{e,\text{up}}=\frac{\bar{x}_e-x_{e}^\circ}{x_{e}^\circ}.
\end{equation*}
The frequency distribution of the optimal adjustments is shown in Fig.~\ref{fig: individual adjustment}. For a line $e$, positive adjustments depict $\gamma_{e,\text{up}}$ values while $\gamma_{e,\text{down}}$ corresponds to negative adjustments. The expectation that increasing $S'$ would require a broader range of control capabilities carries through for $S'$ up to 2000. However, another surprising observation is made for $S'=5000$: the solution of (P2) almost abandons reactance-increase capabilities, and rather attains flow controllability by optimally reducing the reactance of individual lines over wider ranges. Analyzing the shift in solution in the extreme case of $S'=5000$ motivates further research. 

%%%%%%%%%%%%%%%%%%%%%%%%%%%%%%%%%%%%%%%%
\subsection{AC PF based evaluation of power flow control capabilities}
Our empirical tests under DC PF modeling strengthened the analytical findings by demonstrating that desired flow patterns can be realized on a power network by varying the reactances for a subset of lines. Recall that the dataset of flow pattern was generated by solving (P1), which is a relaxation of DC OPF, that does not account for losses or reactive power demands and flows. The flows themselves are therefore infeasible under AC PF settings~\cite{Kyri21DCinfeasible}. Nevertheless, in this section, we evaluate how close the AC power flows can be brought to the flow scenarios $\bef_s$ by varying line reactances. To contextualize the results we solved DC and AC PF for the nominal loads, generation, and network parameters in the test system and found that the average of line flow discrepancies $|f_e^{\mathrm{AC}}-f_e^{\mathrm{DC}}|/\bar{f}_e$ was approximately $1.14\%$. So a mismatch in the AC power flows, after adjusting the line reactances, and desired setpoints $\bef_s$ in that numerical range shall be considered a perfect realization.

We abstractly represent the sending-end active power flow vector using the AC power flow mapping as $\bef_s^{\mrAC}(\bx)=\mathrm{ACPF}(\bp^\mrd_s,\bq^\mrd_s,\bp^\mrg_s;\bx)$, where the power demands and generator active power setpoints are selected from our dataset, and generator voltage setpoints are set to $1$~p.u., and one generator bus is set as the slack bus. With the nominal line reactances, solving AC PF for all $S$ scenarios yielded a flow mismatch $|f_{s,e}^{\mathrm{AC}}(\bx^\circ)-f_{s,e}|/\bar{f}_e$, averaged over $E$ lines and $S$ scenarios, as $8.81\%$. We then investigated the extent to which the mismatch can be reduced if the reactances for all $46$ lines could be controlled. Setting $K=46$ in (P2) provides us, among other things, the optimal reactance vector $\bx_s^{\mrDC}$ for each scenario, where we have added the superscript to emphasize that $\bx_s^{\mrDC}$ is deemed optimal under DC PF settings. We carried out the ensuing analysis for $50$ scenarios that were chosen to represent a varying level of difficulty. The selected instances correspond to $s=1, 101,\dots,4901.$ Recall that, if one sets the line reactances to $\bx_s^{\mrDC}$, the power injections $(\bp^\mrg_s-\bp^\mrd_s)$ would result in the line flows being precisely equal to $\bef_s$ under DC PF. Under AC PF, however, the average mismatch $|f_{s,e}^{\mathrm{AC}}(\bx^\mrDC_s)-f_{s,e}|/\bar{f}_e$, while dropping from $8.81\%$ at nominal reactance, remains at $4.21\%$. The reduction in mismatch is appreciable since the reactance adjustment was calculated by (P2) without accounting for AC PF physics.

For evaluating the added improvement in power flow control by adjusting line reactances explicitly under AC PF considerations, we solved the following optimization routine independently for selected scenarios
\begin{subequations}\label{eq:P4}
\begin{align}
\min_{\bx} \|\bef_s^\mrAC(\bx)-\bef_s\|_2^2+&w\|\bx-\bx_s^{\mrDC}\|_2^2\qquad\textrm{(P3)}\notag\\
\textrm{s.to}\quad \bef_s^{\mrAC}(\bx)&=\mathrm{ACPF}(\bp^\mrd_s,\bq^\mrd_s,\bp^\mrg_s;\bx),\label{seq:P3a}\\
\bx&\geq \mathbf{0},\label{seq:P4d}  
\end{align}
\end{subequations}
where the weight $w$ is used to balance the regularization term that injects useful prior information $\bx^\mrDC_s$ available from solving (P2). Despite the nonconvexity of (P3), it can be efficiently solved by developing gradient-based methods such as the ones reported in~\cite{Babak25tpwrs-flows}. Since numerical tractability was not our emphasis, we solved (P3) for the selected 50 scenarios using the nonlinear least-squares minimization function $\mathrm{lsqnonlin}$ from MATLAB Optimization Toolbox. Since the function can only handle linear constraints, we augmented the cost function to model the AC PF bus power balance mismatches, thereby relaxing~\eqref{seq:P3a}. The AC PF mismatch penalty term was weighted such that the active and reactive power mismatch at all buses was less than $10^{-5}$~p.u. Note that as $w\rightarrow\infty$, the minimizer of (P3) denoted as $\bx^\star_s(w)$ would become equal to $\bx_s^\mrDC$. Small values of $w$, on the other hand, allow aggressive changes in line reactances to minimize the AC line power flow mismatch.
\begin{figure}[h]
    \centering
    \includegraphics[width=0.9\linewidth]{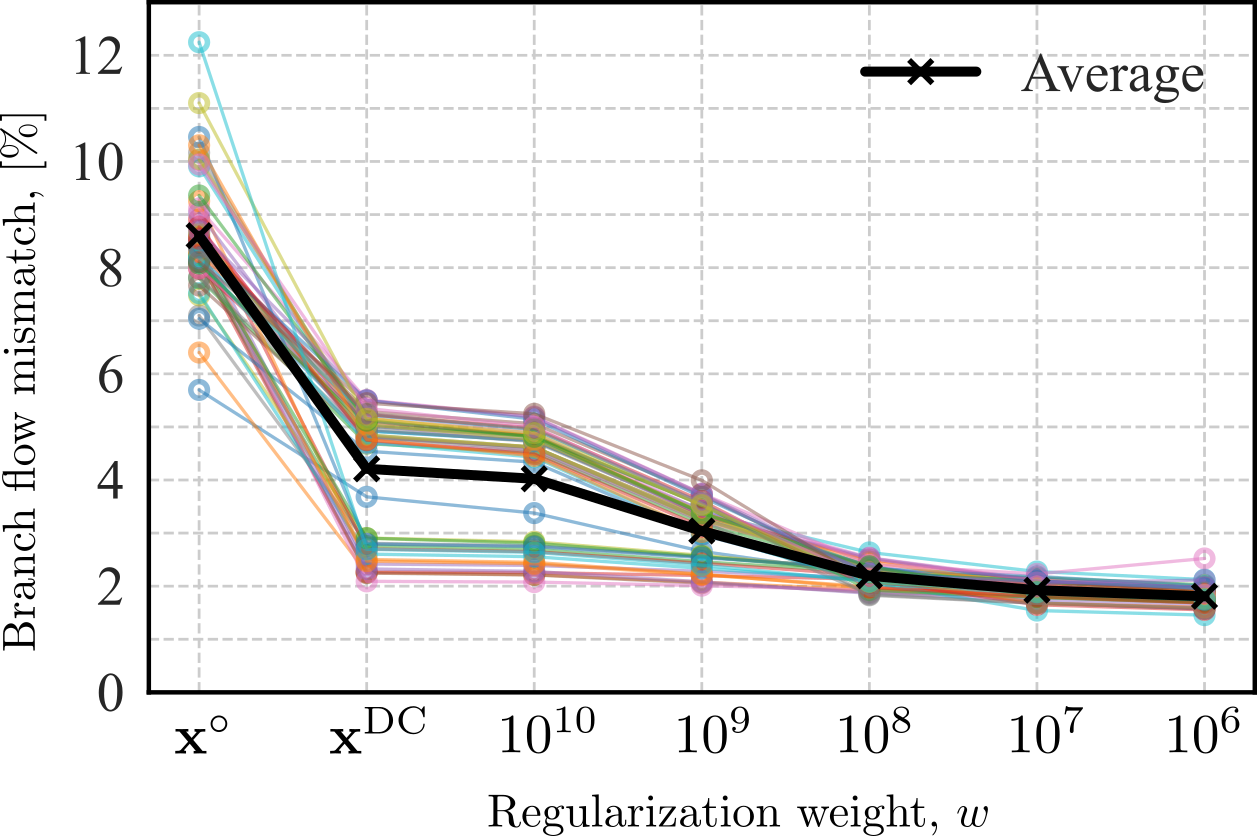}
    \caption{Average (over lines) mismatch between desired line flow values $\bef_s$ and realized line flows under AC PF for varying values of line reactance. Mismatches at nominal reactance, DC setpoint from (P2), and setpoints from (P3) for varying $w$ are plotted. Each trendline represents one scenario, while the bold line provides the average over all 50 scenarios. }
    \label{fig:AC_experiment}
\end{figure}

Figure~\ref{fig:AC_experiment} shows the percentage mismatch in line flows averaged over the 46 lines and 50 scenarios, when the AC power flow $\bef_s^\mrAC$ was computed for the nominal line reactance $\bx^\circ$, reactance setpoint from (P2) $\bx_s^\mrDC$, and the reactance values obtained from (P3) for different values of $w$. We observed that by reducing $w$, and hence allowing (P3) to deviate from the reactances suggested by the DC setup, the average branch flow mismatch could be lowered to $1.8\%$, which can be deemed as a near-perfect match. 

\section{Conclusion and Future Work}
This work has put forth a comprehensive analytical and empirical investigation into the system-wide impact of varying transmission line reactances towards power flow control. Under simplified DC PF settings, it has been shown that arbitrary flow patterns that do not constitute circulating power flows can be realized on a power network if all line reactances are controllable, even if restricted to be inductive. Empirical tests with the DC PF model further reveal that a significant level of flow controllability is achievable with controllers placed on fewer lines. Graph-theoretic insights are obtained on determining which lines would not benefit from flow controllers. Numerical tests under AC PF settings corroborate the validity of our analytical and empirical findings obtained using the simplified DC PF model, and demonstrate the capability of obtaining desired line flow patterns while accounting for detailed AC power flow. This work lays down the algorithmic framework for enabling future research on optimal siting, sizing, and operation of power flow controllers while reporting several observations that warrant further investigation.

\bibliography{myabrv,GETreferences}
\bibliographystyle{IEEEtran}

\end{document}

%% file: StylesMKS.tex
\DeclareMathOperator{\rank}{rank}
\DeclareMathOperator{\range}{range}

\DeclareMathOperator{\nullspace}{null}

\DeclareMathOperator{\diag}{dg}

\newtheorem{proposition}{Proposition}

\newtheorem{property}{Property}

\newcommand \bzero{\mathbf{0}}
\newcommand \bone{\mathbf{1}}

\newcommand \bd{\mathbf{d}}
\newcommand \be{\mathbf{e}}
\newcommand \bef{\mathbf{f}} % exception since \bf is taken

\newcommand \bn{\mathbf{n}}

\newcommand \bp{\mathbf{p}}
\newcommand \bq{\mathbf{q}}

\newcommand \bu{\mathbf{u}}

\newcommand \bx{\mathbf{x}}

\newcommand \bA{\mathbf{A}}

\newcommand \bN{\mathbf{N}}

\newcommand \bQ{\mathbf{Q}}

\newcommand \btheta{\boldsymbol{\theta}}

\newcommand \mrd{\mathrm{d}}

\newcommand \mrg{\mathrm{g}}

\newcommand \mrp{\mathrm{p}}
\newcommand \mrq{\mathrm{q}}

\newcommand \mrAC{\mathrm{AC}}
\newcommand \mrDC{\mathrm{DC}}

\newcommand \mcC{\mathcal{C}}

\newcommand \mcE{\mathcal{E}}

\newcommand \mcG{\mathcal{G}}

\newcommand \mcN{\mathcal{N}}

%% file: Main.bbl
\begin{thebibliography}{10}
\providecommand{\url}[1]{#1}
\csname url@samestyle\endcsname
\providecommand{\newblock}{\relax}
\providecommand{\bibinfo}[2]{#2}
\providecommand{\BIBentrySTDinterwordspacing}{\spaceskip=0pt\relax}
\providecommand{\BIBentryALTinterwordstretchfactor}{4}
\providecommand{\BIBentryALTinterwordspacing}{\spaceskip=\fontdimen2\font plus
\BIBentryALTinterwordstretchfactor\fontdimen3\font minus \fontdimen4\font\relax}
\providecommand{\BIBforeignlanguage}[2]{{%
\expandafter\ifx\csname l@#1\endcsname\relax
\typeout{** WARNING: IEEEtran.bst: No hyphenation pattern has been}%
\typeout{** loaded for the language `#1'. Using the pattern for}%
\typeout{** the default language instead.}%
\else
\language=\csname l@#1\endcsname
\fi
#2}}
\providecommand{\BIBdecl}{\relax}
\BIBdecl

\bibitem{DOE2022GET}
{U.S. Department of Energy}, ``Grid enhancing technologies: A case study on ratepayer impact,'' U.S. Department of Energy, Washington, DC, Tech. Rep., 2022.

\bibitem{Ardakani24survey}
M.~S.-A. Omid~Mirzapour, Xinyang~Rui, ``Grid-enhancing technologies: Progress, challenges, and future research directions,'' \emph{Electric Power Systems Research}, vol. 230, p. 110304, 2024.

\bibitem{Gentle2024PFC}
J.~Gentle, S.~M.~S. Alam, M.~Sun, Z.~Priest, A.~D. Rosso, D.~Schweer, and S.~Guggilam, ``Implementation and operation of power flow control solutions for transmission systems,'' Idaho National Laboratory and EPRI, Technical Report INL/RPT-24-78148, 2024.

\bibitem{GaN3kV}
Y.~Guo, Y.~Qin, M.~Xiao, M.~Porter, Q.~Song, D.~Popa, L.~Efthymiou, K.~Cheng, I.~Kravchenko, L.~Shao, H.~Wang, F.~Udrea, and Y.~Zhang, ``Enhancement-mode gan monolithic bidirectional switch with breakdown voltage over 3.3 kv,'' \emph{{IEEE} Electron Device Lett.}, vol.~46, no.~4, pp. 556--559, 2025.

\bibitem{JahnsGaN}
J.~Liu, H.~Kim, A.~J. Young, Y.~Jiao, E.~Persson, Z.~Xiao, B.~Pandya, T.~M. Jahns, and M.~Imam, ``Gan bidirectional switches: Device technology, applications, and future prospects,'' \emph{{IEEE} Trans. Power Electron.}, pp. 1--14, 2025.

\bibitem{PowerRoute}
R.~P. Kandula, A.~Iyer, R.~Moghe, J.~E. Hernandez, and D.~Divan, ``Power router for meshed systems based on a fractionally rated back-to-back converter,'' \emph{{IEEE} Trans. Power Electron.}, vol.~29, pp. 5172--5180, 2014.

\bibitem{ACFacts}
J.~C. Rosas-Caro, J.~M. Ramirez, and F.~Z. Peng, ``Simple topologies for ac-link flexible ac transmission systems,'' in \emph{IEEE Bucharest PowerTech}, 2009, pp. 1--8.

\bibitem{PowerSeControl}
F.~Mancilla-David and G.~Venkataramanan, ``Modeling and control of the static synchronous series compensator under different operating modes,'' in \emph{IEEE Power Electronics Specialists Conference}, 2007, pp. 2443--2449.

\bibitem{EvalSeries}
F.~Mancilla-David, S.~Bhattacharya, and G.~Venkataramanan, ``A comparative evaluation of series power-flow controllers using dc- and ac-link converters,'' \emph{{IEEE} Trans. Power Delivery}, vol.~23, pp. 985--996, 2008.

\bibitem{Tongxin16tpwrs-flexibility}
J.~Zhao, T.~Zheng, and E.~Litvinov, ``A unified framework for defining and measuring flexibility in power system,'' \emph{{IEEE} Trans. Power Syst.}, vol.~31, pp. 339--347, 2016.

\bibitem{Liu2018gridflexibility}
J.~Li, F.~Liu, Z.~Li, C.~Shao, and X.~Liu, ``Grid-side flexibility of power systems in integrating large-scale renewable generations: A critical review on concepts, formulations and solution approaches,'' \emph{Renewable and Sustainable Energy Reviews}, vol.~93, pp. 272--284, 2018.

\bibitem{Tomsovic17tpwrs}
X.~Zhang, K.~Tomsovic, and A.~Dimitrovski, ``Security constrained multi-stage transmission expansion planning considering a continuously variable series reactor,'' \emph{{IEEE} Trans. Power Syst.}, vol.~32, pp. 4442--4450, 2017.

\bibitem{ChananSingh2001}
Y.~Ou and C.~Singh, ``Improvement of total transfer capability using {TCSC} and {SVC},'' in \emph{Power Engineering Society Summer Meeting. Conference Proceedings}, vol.~2, 2001, pp. 944--948.

\bibitem{ChananSingh2002tpwrs}
------, ``Assessment of available transfer capability and margins,'' \emph{{IEEE} Trans. Power Syst.}, vol.~17, pp. 463--468, 2002.

\bibitem{MATPOWER}
R.~D. Zimmerman, C.~E. Murillo-Sanchez, and R.~J. Thomas, ``{MATPOWER}: steady-state operations, planning and analysis tools for power systems research and education,'' \emph{{IEEE} Trans. Power Syst.}, vol.~26, pp. 12--19, 2011.

\bibitem{Kyri21DCinfeasible}
K.~Baker, ``Solutions of {DC} {OPF} are never {AC} feasible,'' in \emph{Proc. ACM International Conference on Future Energy Systems}, ser. e-Energy '21.\hskip 1em plus 0.5em minus 0.4em\relax New York, NY, USA: Association for Computing Machinery, 2021, p. 264–268.

\bibitem{Babak25tpwrs-flows}
B.~Taheri and D.~K. Molzahn, ``Ac-informed dc optimal transmission switching problems via parameter optimization,'' \emph{{IEEE} Trans. Power Syst.}, vol.~40, no.~6, pp. 5422--5433, 2025.

\end{thebibliography}
